\documentclass[12pt,letterpaper]{amsart}
\usepackage{euler, epic,eepic,latexsym, amssymb, amscd, amsfonts, xypic, floatflt}
\input xy
\xyoption{all}

 \newlength{\baseunit}               
 \newcount{\numlines}                
\newcommand{\point}{\vspace{2mm}\par \noindent
  \refstepcounter{subsection}{\thesubsection.} } 
\newcommand{\tpoint}[1]{\vspace{2mm}\par \noindent \refstepcounter{subsection}{\thesubsection.} 
  {\bf #1. ---} }
\newcommand{\epoint}[1]{\vspace{2mm}\par \noindent \refstepcounter{subsection}{\thesubsection.} 
  {\em #1.} }
\newcommand{\bpoint}[1]{\vspace{2mm}\par \noindent \refstepcounter{subsection}{\thesubsection.} 
  {\bf #1.} }

\newcommand{\Z}{\mathbb{Z}}
\newcommand{\A}{\mathbb{A}}

\newcommand{\Q}{\mathbb{Q}}
\newcommand{\Pbar}{\overline{P}}
\newcommand{\Qbar}{\overline{Q}}

\newcommand{\C}{\mathbb{C}}
\newcommand{\F}{\mathbb{F}}
\newcommand{\Fbar}{\overline{\F}}

\newcommand{\proj}{\mathbb P}

\newcommand{\Spec}{\operatorname{Spec}}

\newcommand{\lremind}[1]{{\bf[label:  #1]}}
\newcommand{\notation}[1]{}
\renewcommand{\lremind}[1]{{}}

\newcommand{\cut}[1]{}

\begin{document}
\pagestyle{plain}
\title{{\large {Mn\"ev-Sturmfels universality for schemes}}
}
\author{Seok Hyeong Lee and Ravi Vakil}\thanks{The first author was supported by Samsung Scholarship. The second author was supported by NSF grant DMS-1100771.}
\address{Department of Mathematics, Stanford University, Stanford CA~94305--2125}
\email{lshyeong@stanford.edu, vakil@math.stanford.edu}
\date{July 10, 2012.}
\begin{abstract}
We prove a scheme-theoretic version of Mn\"ev-Sturmfels Universality, suitable to be used in the proof of Murphy's Law in Algebraic Geometry \cite[Main Thm.~1.1]{v}.  Somewhat more precisely, we show that any singularity type of finite type over $\Z$ appears on some incidence scheme of points and lines, subject to some particular further constraints.

This paper is dedicated to Joe Harris on the occasion of his birthday, with warmth and gratitude.
\end{abstract}
\maketitle
\tableofcontents

{\parskip=12pt 

\section{Introduction}

\point \lremind{d}\label{d}Define an equivalence relation $\sim$ on pointed schemes generated by the following: if $(X,P) \rightarrow (Y,Q)$ is a smooth morphism of pointed schemes ($P \in X$, $Q \in Y$) --- i.e.\ a smooth morphism $\pi: X \rightarrow Y$ with $\pi(P) = Q$ ---  then $(X,P) \sim (Y,Q)$.  We call equivalence classes  {\em singularity types}, and we call pointed schemes {\em singularities}.  We say that {\em Murphy's Law} holds for a (moduli) scheme $M$ if every singularity type appearing on a finite type scheme over $\Z$ also appears on $M$.  (This use of the phrase ``Murphy's Law'' is from \cite[\S 1]{v}, and earlier appeared informally in \cite[p.~18]{hm}. Folklore ascribes it to Mumford.)

\epoint{Definition}  \label{d:is}Define an {\em incidence scheme} of points and lines in $\proj^2_{\Z}$ as a locally closed subscheme of $(\proj^2_{\Z})^M \times (\proj^{2 \vee}_{\Z})^N = \{ p_1, \dots, p_M, l_1, \dots, l_N \}$ parametrizing $M$ labeled points and $N$ labeled lines, satisfying the following conditions.
\begin{enumerate}
\item[(i)] $p_1 = [0,0,1]$, $p_2=[0,1,0]$, $p_3 = [1,0,0]$, $p_4 = [1,1,1]$.
\item[(ii)] We are given some specified incidences:  for each pair $(p_i, l_j)$, either $p_i$ is required to lie on $l_j$, or $p_i$ is required not to lie on $l_j$.
\item[(iii)] The marked points are required to be distinct, and the marked lines are required to be distinct.
\item[(iv)] Given any two marked lines, there is a marked point required to be on both of them (necessarily unique, given (iii)).
\item[(v)] Each marked line contains at least three marked points.
\end{enumerate}

Note that even though our definition over $\Z$, these conditions may
force us into positive characteristic.  For instance, the Fano plane
point-line configuration would force us into characteristic $2$.

The goal of this paper is to establish the following.

\tpoint{Mn\"ev-Sturmfels Universality Theorem for Schemes}  {\em The disjoint union of all incidence schemes (over all possible $M$, $N$, and data of the form described in Definition~\ref{d:is}) satisfies Murphy's Law.}\label{t:main}\lremind{t:main}


Theorem~\ref{t:main} appeared as \cite[Thm.~3.1]{v}, as an essential
step in proving \cite[Main Thm.~1.1]{v}, which stated that many
important moduli spaces satisfy Murphy's Law.  A number of readers of
\cite{v} have pointed out to the second author that the references
given in \cite{v} and elsewhere in the literature do not establish the
precise statement of Theorem~\ref{t:main}, and that it is not clear
how to execute the glib parenthetical assertion (``The only subtlety
...,'' \cite[p.~577, l.~3-5]{v}) to extend Lafforgue's argument
\cite[Thm.~I.14]{l} to obtain the desired result.  (In particular,
our problem is that Lafforgue does not require (iv) in his moduli
spaces, as it is not needed for his purposes.  If we then add
marked points to all pairwise intersections of lines in Lafforgue's construction, then 
it is not clear that (iii) holds in the configurations he constructs,
and we suspect it does not always hold.)

This paper was written in order
to fill a possible gap in \cite{v}, or at least to clarify details of
an important construction.  
We hope this paper will be of use to those studying the singularities
of moduli spaces not covered by \cite{v} (the moduli space of vector
bundles, \cite{p}, or the Hilbert scheme of points, \cite{e}, say).
Although no one familiar with this area
would doubt that Theorem~\ref{t:main} holds, or how the general idea
should go, we will see that some
care is needed to rigorously establish it.  In particular, our
argument is characteristic-dependent.

\epoint{Key features in the argument} 
Given any polynomials $f_1, \dots, f_r \in   \Z[x_1, \dots, x_n]$, our
goal is to build a smooth cover of
\[
\Spec \Z[x_1, \dots, x_n]/(f_1, \dots, f_r)
\]
by (open subsets of) incidence schemes, by encoding the variables and relations in incidence relations. We
build the relations by combining ``atomic'' calculations encoding
equality, negation, addition, and multiplication. We point out new
features of the argument we use, in order to ensure \ref{d:is}(iii) in
particular. We perform each ``atomic'' calculation on a separate line
of the plane, to avoid having too many important points on a single
line, because points on a line  must be shown to not overlap.  We need various cases to deal with when
the ``variable'' in question is ``near'' $0$ or $1$ (i.e.\ has value
$0$ or $1$ at a given geometric point $Q$ of $\Z[x_1, \dots, x_n]/(f_1,
\dots, f_r)$, but is not required to have that value ``near''
$Q$). Furthermore, the ``usual'' construction of addition and
multiplication runs into problems in characteristic $2$ due to
unintended coincidences of points, so some care is required in this
case (see \S \ref{s:char2}).

\epoint{Algebro-geometric history} Vershik's ``universality''
philosophy (e.g.\ \cite[Sect.~7]{vershik}) has led to a number of
important constructions in many parts of mathematics. One of the most
famous is Mn\"ev's Universality Theorem~\cite{m1, m2}.  It was
independently proved by Bokowski and Sturmfels~\cite{bs, s1, s2}.  We
follow Belkale and Brosnan \cite[\S 10]{bb} in naming the result after
both Mn\"ev and Sturmfels.  (The idea is more ancient; von Staudt's
``algebra of throws'' goes back at least to \cite{mac}, see also
\cite{k}.)

Lafforgue outlined a proof of a scheme-theoretic version in
\cite[Thm.~I.14]{l}. Keel and Tevelev used this construction in
\cite{kt} (see \S 1.8 and Theorem~3.13 of that article). Another
algebro-geometric application of Mn\"ev's theorem (this time in its
manifestation in the representation problem of matroids) was Belkale
and Brosnan's surprising counterexample to a conjecture of Kontsevich,
\cite{bb}.  More recently Payne applied this construction in \cite{p} to toric vector
bundles, and Erman applied it in \cite{e} to the Hilbert scheme of
points.  (These examples are representative but not exhaustive.)

\noindent {\bf Acknowledgments.}
We thank the two referees for their  thoughtful suggestions.  The second author thanks A. J. de Jong, M. Roth, and of course
J. Harris.  

\section{Structure of the construction}

\bpoint{Strategy} \label{strategy} Fix a singularity $(Y, Q)$ of
finite type over $\Spec \Z$. We will show that there exists a point
$P$ of some incidence scheme $X$ (i.e.\ some configuration of points
and lines, as described in Definition~\ref{d:is}), along with a smooth
morphism  $\pi: (X,P) \rightarrow (Y,Q)$  of pointed schemes. Because smooth morphisms are open, it suffices to deal with the case where $Q$ is a {\em closed} point of $Y$. Then the residue field $\kappa(Q)$ has finite characteristic $p$. (The reduction to characteristic $p$ is not important; it is done to allow us to construct a configuration over a fixed infinite field. Those interested only in the characteristic $0$ version of this result will readily figure out how to replace $\Fbar_p$ with $\Q$ or $\C$.)

By replacing $Y$ with an affine neighborhood of $Q$, we may assume $Y$ is affine, say $Y = \Spec \Z[x_1, \dots, x_n] / (f_1, \dots, f_r)$. The morphism $$\Spec \Fbar_p[x_1, \dots, x_n] / (f_1, \dots, f_r) \rightarrow \Spec \F_p[x_1, \dots, x_n] / (f_1, \dots, f_r)$$ is surjective by the Lying Over Theorem. Choose a pre-image $\Qbar \in \Spec \Fbar_p[x_1, \dots, x_n] / (f_1, \dots, f_r)$ of $Q$ --- say the (closed) point $(x_1, \dots, x_n) = (q_1, \dots, q_n)$, where $\vec{q} \in \Fbar_p^n$.

We make the following constructions.

\noindent (a)
We describe  a configuration of points and lines over $\Fbar_p$, which is thus an $\Fbar_p$-valued point $\Pbar$ of an incidence scheme $X$.

\noindent (b) The incidence scheme $X$ will be an open subscheme of an
affine scheme $X'$, and we construct (a finite number of) coordinates
on $X'$, which we name $X_1$, \dots, $X_n$, $Y_1$, \dots, $Y_s$
subject only to the relations $f_i(X_1, \dots, X_n)=0$ ($1 \leq i \leq
r$).  We thus have a smooth morphism $\pi: X \rightarrow Y$ given by
$X_i \mapsto x_i$.  Letting $X_K = X \times_{\Spec \Z} \Spec K$ for $K=\F_p$ and $\Fbar_p$, and similarly for $Y_K$ and $\pi_K$, we have a diagram:
\[
 \xymatrix{
  X_{\Fbar_p} \ar[r] \ar[d]_{\pi_{\Fbar_p}} & X_{\F_p} \ar[r]
  \ar[d]_{\pi_{\F_p}} & X \ar@{^(->}[r]^{\text{open}} \ar[d]_{\pi} &
  X' \ar@{=}[r] & Y \times_{\Z} \A^s_{\Z} \ar[dll]
  \\
  Y_{\Fbar_p} \ar[r] & Y_{\F_p}\ar[r] & Y}
\]
In the course of the construction, we will not explicitly name the
variables $Y_j$, but whenever a free choice is made this corresponds to
adding a new variable $Y_j$.

\noindent (c) We will have $\pi_{\Fbar_p}( \Pbar) = \Qbar$.  Thus the image of $\pi$ includes $Q$.

\bpoint{Notation and variables for the incidence scheme}

The traditional (and only reasonable) approach is to construct a configuration of points and lines encoding this singularity, by encoding the ``atomic'' operations of equality, negation, addition, and multiplication.  The most difficult desideratum is \ref{d:is}(iii).

Our incidence scheme will parametrize points and lines of the following form. In the course of this description we give names to the relevant types of points and lines, and give our chosen coordinates. We will later describe our particular point $\Pbar \in X_{\Fbar_p}$.

The first type of points are $p_1$ through $p_4$ (see
\ref{d:is}(i)). We call these {\em anchor points}. We interpret
$\proj^2$ in the usual way:  $p_2 p_3$ is the line at infinity, and $p_1$ is the origin; lines through $p_3$ are called horizontal.  The first type of lines in our incidence scheme are the lines $p_i p_j$.  We call these {\em anchor lines}.

The next type of line, which we call {\em variable-bearing lines},
will be required to pass through $p_3=[1,0,0]$ (they are
``horizontal''), and not through $p_1$,$p_2$, or $p_4$. Each
variable-bearing line is parametrized by where it meets the $y$-axis
(which is finite, as the lines do not pass through
$p_2=[0,1,0]$). Thus for each variable-bearing line $l_i$, we have a
coordinate $y_i$.  In order to satisfy \ref{d:is}(iii), we will always
arrange that the $y_i$ are distinct and not $0$ or $1$.  These
$y_i$ will be among the $Y_j$ of \ref{strategy}(b) above.

Each variable-bearing line $l_i$ has a {\em framing-type} $Fr_i$,
which is a size two subset of 
$\{ -1, 0, 1\}$ if $p>2$, and
of $\{ 0, 1, j \}$ where $j$ is a chosen solution of $j^2+j-1=0$ (see
\S \ref{s:if2} for more)
 if $p=2$.  Each
variable-bearing line $l_i$ contains (in addition to $p_3$) the
following three distinct marked points: \begin{itemize}
\item two {\em framing points} $P_{i,s}$, where $s \in Fr_i$; and 
\item one {\em variable-bearing point} $V_i$.
\end{itemize}
The point $P_{i,s}$ we parametrize by its $x$-coordinate, which we
confusingly name $y_{i,s}$ (because it will be one of the ``free'' variables $Y_i$
of \ref{strategy}(b)). The variable-bearing point $V_i$ we parametrize using
the isomorphism $l_i \rightarrow \proj^1$ obtained by sending $p_3$ to
$\infty$, and $P_{i,s}$ to $s$ for $s \in Fr_i$. We denote this coordinate $x_i$. 
(In our construction, these coordinates will be  either  among those
$X_i$ of \ref{strategy}(b) above, or will be determined by the other
variables.)  A variable-bearing line over $\Fbar_p$ of framing-type $Fr_i$, whose
variable-bearing point carries the variable $x_i = q \in \Fbar_p$, we
will call a $( Fr_i, q)$-line or a $( Fr_i, x_i)$-line.  

We have a number of additional configurations of points and lines,
called {\em connecting configurations}, which are required to contain
a specified subset of the above-named points, and required to {\em
  not} contain the rest. These will add additional free variables
(which, in keeping with \ref{strategy}(b) above, we call $Y_j$ for an
appropriate $j$), and will (scheme-theoretically)  impose a single constraint upon the $x$-variables: 
\begin{itemize}
\item $x_a = x_{b}$ (an {\em equality configuration})
\item $x_a = -x_{b}$ (a {\em negation configuration}),
\item $x_{a} + x_{b} = x_{c}$ (an {\em addition configuration}), or
\item $x_{a}  x_{b} =  x_{c}$ (a {\em multiplication configuration}). 
\end{itemize}

Finally, for each pair of above-named lines that do not have an
above-named point contained in both, we have an additional marked
point at their intersection (in order that \ref{d:is}(iv) holds),
which we call {\em bystander points}. In our
construction, we never have more than two lines meeting at a point
except at the previously-named points, and the only lines passing through
the previously-named points are the ones specified above.  The name
``bystander points'' reflects the fact that they play no further
role, and no additional variables are needed to parametrize them.

\bpoint{Reduction to four problems}

\label{s:four}We reduce Theorem~\ref{t:main} to four atomic  problems.

We construct the expression for each $f_i$ sequentially, starting with
the variables $x_i$ and the constant $1$ (which we name $x_0$ to
simplify notation later), and using negation of one term
or addition or multiplication of two terms at each step. Somewhat more
precisely, we make a finite sequence of intermediate expressions, where each
expression is the negation, sum, or product of earlier (one or two)
expressions. We assign new variables $x_{n+1}, x_{n+2}, \dots$ for
each new intermediate expression. Those additional variables $x_k$ will come
along  with a single equation ---  negation, addition, or
multiplication --- describing how $x_k$ is obtained from its
predecessor(s). In case of sum and product, we additionally require two predecessors to be different, even in case of adding or multiplying same expression. Finally, for the variable $x_a$ representing the
final expression $f_i$ (one for each $f_i$), we add the equation $x_a = 0$. These
simple equations (which we call $g_i$) are equivalent to our original equations
$f_i=0$, so we have (canonically)
\begin{equation}\label{e:fg}
 \Z[x_1, \dots, x_n]/(f_1, \dots, f_r) \cong \Z[x_1, \dots, x_n,
 x_{n+1}, \dots, x_m]/(g_1, g_2, \dots, g_{r'})\end{equation}
where each $g_j$ is of the form $x_a - x_b$, $x_a + x_b$,  $x_a + x_b -
x_c$, $x_a  x_b -  x_c$, or $x_a$.   As one example of this procedure:
$$
\frac { \Z [x_1, x_2, x_3]} { (x_1 x_2 +x_3^2-2)}
\cong 
\frac { \Z [x_1, x_2, x_3, x_4, x_5, x_6, x_7, x_8, x_9, x_{10}]}
 { (x_4-x_1 x_2, x_5-x_3, x_6-x_3x_5, x_7-x_4-x_6, x_8-1, x_9-1,
   x_{10}-x_8-x_9, x_7 -x_{10})}.
$$

We now construct our configuration over $\Fbar_p$.  Via \eqref{e:fg},
we interpret $\Qbar$ as a geometric point of $\Spec \Z[x_1, \dots,
x_m]/(g_1, g_2, \dots, g_{r'})$, and we let $q_i \in \Fbar_p$ be
coordinates of $x_i$ for all $i$. 
For each $i \in \{1, \dots, m \}$, we choose two distinct $q$'s in $\{-1, 0, 1\}$ or $\{0, 1, j \}$ (according to whether $p > 2$
or $p=2$) distinct from $q_i$; this will be the framing-type $Fr_i$ of $x_i$.
We place (generally chosen) variable-bearing
lines $l_i$, one for each $x_i$ for $i \in \{1, \dots, m\}$, with
framing points (corresponding to the framing-type $Fr_i$) chosen generally on $l_i$, then with
variable-bearing points chosen so that the coordinate of the variable
point for the line $l_i$ is $q_i$. 

We then sequentially do the following for each simple equation
$g_j$. For each $g_j$ involving variables $x_a, x_b$ (and possibly
$x_c$), we place a corresponding configuration joining variable-bearing points
for those variables and enforcing (scheme-theoretically) the equation
$g_j$.  We will do this in such a way that the connecting
configuration passes through no points or lines it is not supposed
to. We will of course do this by a general position argument.

We are thus reduced to four problems, which we describe below (with
italicized titles) after setting the stage for them.  Suppose we are
given a configuration of points and lines in the plane, including the
anchor points $p_i$, and the anchor lines $p_i p_j$ ($1 \leq i < j
\leq 4$) (and hence implicitly a point of some incidence scheme). Note
that this incidence scheme is quasiaffine, say $U \subset \Spec A$:
\begin{itemize}
\item  the
non-vertical lines (those non-anchor lines not containing $p_2 = [
0,1,0]$) $y=mx+b$ are parametrized by $m$ and $b$;
\item  vertical lines
(those non-anchor lines passing through $p_2$) $x=a$ are parametrized
by $a$;
\item  those points $(x,y)=(a,b)$ not on the line at infinity are
parametrized by $a$ and $b$;
\item and those non-anchor points $[1,c,0]$ on
the line at $\infty$ are parametrized by $c$. \end{itemize} The conditions of
\S \ref{d:is} are clearly locally closed.

Here now are the four problems.

{\em Equality problem.}  If we have two variable-bearing lines $l_a$
and $l_b$ with coordinates $x_a$ and $x_b$, we must show that we may
superimpose an equality configuration (i.e.\ add more points and lines),
where except for the framing and variable-bearing points on these
two lines $l_a$ and $l_b$, no point of the additional configuration lies on any
pre-existing lines, and no line in the additional configuration passes
through any pre-existing points (including pre-existing bystander
points --- pairwise intersections of pre-existing lines). Furthermore,
the addition of this configuration must add only open conditions for
added free variables and $x_a$, $x_b$, and enforce exactly 
(scheme-theoretically) 
the
equation $x_a = x_b$. More precisely, we desire
that the
morphism from the new incidence scheme to the old one is of the
following form:
\[
\xymatrix{
U' \ar[d] \ar@{^(->}[r]^-{\text{open}} & \Spec A[y_1, \dots, y_N] / (x_a - x_b) \ar[d] \\
U \ar[r]  \ar@{^(->}[r]^{\text{open}} & \Spec A
}
\]
(for some value of $N$).

{\em Negation problem:} the same problem, except with $x_a = -x_b$ replacing $x_a = x_b$.

{\em Addition problem:} the analogous problem, except $x_c = x_a + x_b$ ($a \neq b$).

{\em Multiplication problem:} the analogous problem, except $x_c = x_a x_b$ ($a \neq b$).

\section{The configurations}

We now describe the configurations needed to make this work.

\bpoint{Building blocks for the building blocks:  five configurations}

\label{s:atomic}The building blocks we use are shown in
Figures~\ref{f:parallelshift}--\ref{f:genericmultiplication}.  The
figures follow certain conventions.  (See Figure~\ref{f:legend} for a
legend.) Lines that appear horizontal are indeed so --- they are
required to pass through $p_3 = [1,0,0]$.    The
horizontal lines often have (at least) three labeled points, which
suggest an isomorphism with $\proj^1$.  The dashed lines (and marked
points thereon) are those that are in the configuration before we
begin.  The points and lines marked with a box are added next, and
involve free choices (two coordinates for each boxed point, one for the each
boxed horizontal line). The remaining points and lines are then determined.
The triangle indicates the ``goal'' of the construction, if
interpreted as constructing midpoint, addition, multiplication, and so
forth (which is admittedly not our point of view).

\begin{figure}[ht]
\setlength{\unitlength}{0.00083333in}
\begingroup\makeatletter\ifx\SetFigFont\undefined%
\gdef\SetFigFont#1#2#3#4#5{%
  \reset@font\fontsize{#1}{#2pt}%
  \fontfamily{#3}\fontseries{#4}\fontshape{#5}%
  \selectfont}%
\fi\endgroup%
{\renewcommand{\dashlinestretch}{30}
\begin{picture}(3513,1047)(0,-10)
\put(312,936){\blacken\ellipse{60}{60}}
\put(312,936){\ellipse{60}{60}}
\put(87,636){\blacken\ellipse{60}{60}}
\put(87,636){\ellipse{60}{60}}
\put(462,111){\blacken\ellipse{60}{60}}
\put(462,111){\ellipse{60}{60}}
\dashline{60.000}(12,936)(1212,936)
\path(12,711)(162,711)(162,561)
	(12,561)(12,711)
\path(1137,636)(387,636)
\path(612,711)(762,711)(762,561)
	(612,561)(612,711)
\path(386,63)(461,213)(536,63)(386,63)
\put(1437,936){\makebox(0,0)[lb]{{\SetFigFont{8}{9.6}{\rmdefault}{\mddefault}{\updefault}previously constructed lines and points}}}
\put(1437,336){\makebox(0,0)[lb]{{\SetFigFont{8}{9.6}{\rmdefault}{\mddefault}{\updefault}(other points and lines are determined)}}}
\put(1437,636){\makebox(0,0)[lb]{{\SetFigFont{8}{9.6}{\rmdefault}{\mddefault}{\updefault}freely chosen points and horizontal lines}}}
\put(1437,36){\makebox(0,0)[lb]{{\SetFigFont{8}{9.6}{\rmdefault}{\mddefault}{\updefault}``goal''}}}
\end{picture}
}
\caption{Legend for Figures~\ref{f:parallelshift}--\ref{f:genericmultiplication}}\label{f:legend}
\end{figure}
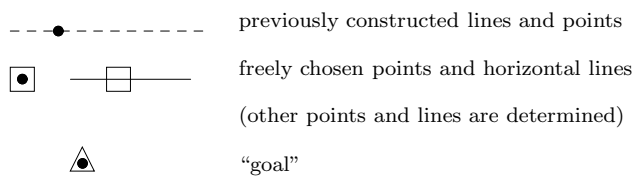

The first building block, {\em parallel shift}, gives the projection
from a point $X$ of
three points $(P_1, V, P_2)$ on the horizontal line $l$ onto $(P_1',
V', P_2')$ on the horizontal line $l'$. Invariance of cross-ratio under projection gives
\[
(P_1, V; P_2, p_3) = (P_1', V'; P_2', p_3)
\]
(where $(\cdot, \cdot; \cdot, \cdot)$ throughout the paper means cross-ratio, or moduli
point in $\mathcal{M}_{0,4}$),
so if $l$ and $l'$ are lines of the same framing-type, with  $(P_1, P_2)$ and
$(P_1', P_2')$ the corresponding framing points and $V$ and $V'$
the  variable-bearing points, then the coordinates
of $V$ and $V'$ are the same (scheme-theoretically). 
Note that we are adding three free variables (two for the point, one
for the line), plus an open condition to ensure no unintended
incidences with preexisting points and lines.

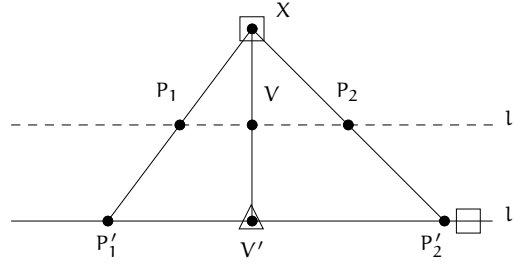
\begin{figure}[ht]
\setlength{\unitlength}{0.00083333in}
\begingroup\makeatletter\ifx\SetFigFont\undefined%
\gdef\SetFigFont#1#2#3#4#5{%
  \reset@font\fontsize{#1}{#2pt}%
  \fontfamily{#3}\fontseries{#4}\fontshape{#5}%
  \selectfont}%
\fi\endgroup%
{\renewcommand{\dashlinestretch}{30}
\begin{picture}(3286,1630)(0,-10)
\put(1512,1456){\blacken\ellipse{60}{60}}
\put(1512,1456){\ellipse{60}{60}}
\put(1062,856){\blacken\ellipse{60}{60}}
\put(1062,856){\ellipse{60}{60}}
\put(612,256){\blacken\ellipse{60}{60}}
\put(612,256){\ellipse{60}{60}}
\put(1512,256){\blacken\ellipse{60}{60}}
\put(1512,256){\ellipse{60}{60}}
\put(1512,856){\blacken\ellipse{60}{60}}
\put(1512,856){\ellipse{60}{60}}
\put(2112,856){\blacken\ellipse{60}{60}}
\put(2112,856){\ellipse{60}{60}}
\put(2712,256){\blacken\ellipse{60}{60}}
\put(2712,256){\ellipse{60}{60}}
\dashline{60.000}(12,856)(3012,856)
\path(1512,256)(1512,1456)
\path(1512,1456)(2712,256)
\path(1512,1456)(612,256)
\path(1437,1531)(1587,1531)(1587,1381)
	(1437,1381)(1437,1531)
\path(2787,331)(2937,331)(2937,181)
	(2787,181)(2787,331)
\path(3012,256)(12,256)
\path(1433,207)(1508,357)(1583,207)(1433,207)
\put(1587,1006){\makebox(0,0)[lb]{{\SetFigFont{8}{9.6}{\rmdefault}{\mddefault}{\updefault}$V$}}}
\put(2037,1006){\makebox(0,0)[lb]{{\SetFigFont{8}{9.6}{\rmdefault}{\mddefault}{\updefault}$P_2$}}}
\put(537,31){\makebox(0,0)[lb]{{\SetFigFont{8}{9.6}{\rmdefault}{\mddefault}{\updefault}$P'_1$}}}
\put(1437,31){\makebox(0,0)[lb]{{\SetFigFont{8}{9.6}{\rmdefault}{\mddefault}{\updefault}$V'$}}}
\put(3087,856){\makebox(0,0)[lb]{{\SetFigFont{8}{9.6}{\rmdefault}{\mddefault}{\updefault}$l$}}}
\put(1662,1531){\makebox(0,0)[lb]{{\SetFigFont{8}{9.6}{\rmdefault}{\mddefault}{\updefault}$X$}}}
\put(912,1006){\makebox(0,0)[lb]{{\SetFigFont{8}{9.6}{\rmdefault}{\mddefault}{\updefault}$P_1$}}}
\put(2562,31){\makebox(0,0)[lb]{{\SetFigFont{8}{9.6}{\rmdefault}{\mddefault}{\updefault}$P'_2$}}}
\put(3087,256){\makebox(0,0)[lb]{{\SetFigFont{8}{9.6}{\rmdefault}{\mddefault}{\updefault}$l'$}}}
\end{picture}
}
\caption{Parallel shift}\label{f:parallelshift}
\end{figure}

The second building block,  {\em midpoint} (Figure~\ref{f:midpoint}), will be used for constructing the midpoint
$M$ of two distinct points $A$ and $B$ on line $l$ (where $p_3$ is
considered as usual to be infinity). This construction will be used
outside characteristic 2.  (In $p=2$, the diagram is misleading:  $XY$
passes through $p_3$, resulting in $M=p_3$.)
We have equality of cross-ratios
\begin{align*}
(A, M; B, p_3) &= (A', M'; B', p_3) \quad \text{(projection
  from $X$)}\\
&= (B, M; A, p_3) \quad \text{(projection from $Y$)}  \\
&= 1/(A, M; B, p_3) \quad \text{(property of cross-ratio)}
\end{align*}
so $(A, M; B, p_3)$ is either $1$ or $-1$. For $p \neq 2$ and $A
\neq B$, it is straightforward to verify that $M \neq p_3$ so $(A, M;
B, p_3) \neq 1$. Thus $(A, M; B, p_3) = -1$, so 
$M$ is the ``midpoint'' of $AB$.   
(More precisely:    given any isomorphism of $l$ with $\proj^1$
identifying $p_3$ with $\infty$, the coordinate of $M$ is the average
of the coordinates of $A$ and $B$.  In classical language, $M$ is the
harmonic conjugate of $p_3$ with respect to $A$ and $B$.)

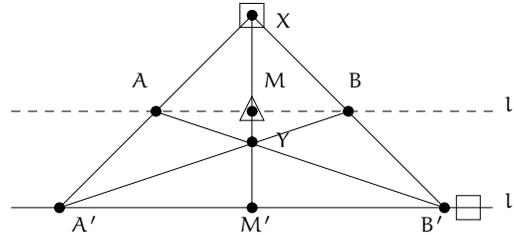
\begin{figure}[ht]
\setlength{\unitlength}{0.00083333in}
\begingroup\makeatletter\ifx\SetFigFont\undefined%
\gdef\SetFigFont#1#2#3#4#5{%
  \reset@font\fontsize{#1}{#2pt}%
  \fontfamily{#3}\fontseries{#4}\fontshape{#5}%
  \selectfont}%
\fi\endgroup%
{\renewcommand{\dashlinestretch}{30}
\begin{picture}(3286,1483)(0,-10)
\put(1512,1381){\blacken\ellipse{60}{60}}
\put(1512,1381){\ellipse{60}{60}}
\put(1512,181){\blacken\ellipse{60}{60}}
\put(1512,181){\ellipse{60}{60}}
\put(1512,781){\blacken\ellipse{60}{60}}
\put(1512,781){\ellipse{60}{60}}
\put(2112,781){\blacken\ellipse{60}{60}}
\put(2112,781){\ellipse{60}{60}}
\put(2712,181){\blacken\ellipse{60}{60}}
\put(2712,181){\ellipse{60}{60}}
\put(312,181){\blacken\ellipse{60}{60}}
\put(312,181){\ellipse{60}{60}}
\put(912,781){\blacken\ellipse{60}{60}}
\put(912,781){\ellipse{60}{60}}
\put(1512,591){\blacken\ellipse{60}{60}}
\put(1512,591){\ellipse{60}{60}}
\path(1512,181)(1512,1381)
\path(1512,1381)(2712,181)
\path(1512,1381)(312,181)
\path(1437,1456)(1587,1456)(1587,1306)
	(1437,1306)(1437,1456)
\path(2787,256)(2937,256)(2937,106)
	(2787,106)(2787,256)
\path(3012,181)(12,181)
\dashline{60.000}(12,781)(3012,781)
\path(912,781)(2712,181)
\path(2112,781)(312,181)
\path(1437,726)(1512,876)(1587,726)(1437,726)
\put(3087,781){\makebox(0,0)[lb]{{\SetFigFont{8}{9.6}{\rmdefault}{\mddefault}{\updefault}$l$}}}
\put(1587,931){\makebox(0,0)[lb]{{\SetFigFont{8}{9.6}{\rmdefault}{\mddefault}{\updefault}$M$}}}
\put(1662,1306){\makebox(0,0)[lb]{{\SetFigFont{8}{9.6}{\rmdefault}{\mddefault}{\updefault}$X$}}}
\put(762,931){\makebox(0,0)[lb]{{\SetFigFont{8}{9.6}{\rmdefault}{\mddefault}{\updefault}$A$}}}
\put(2112,931){\makebox(0,0)[lb]{{\SetFigFont{8}{9.6}{\rmdefault}{\mddefault}{\updefault}$B$}}}
\put(1662,556){\makebox(0,0)[lb]{{\SetFigFont{8}{9.6}{\rmdefault}{\mddefault}{\updefault}$Y$}}}
\put(387,31){\makebox(0,0)[lb]{{\SetFigFont{8}{9.6}{\rmdefault}{\mddefault}{\updefault}$A'$}}}
\put(1437,31){\makebox(0,0)[lb]{{\SetFigFont{8}{9.6}{\rmdefault}{\mddefault}{\updefault}$M'$}}}
\put(2562,31){\makebox(0,0)[lb]{{\SetFigFont{8}{9.6}{\rmdefault}{\mddefault}{\updefault}$B'$}}}
\put(3087,181){\makebox(0,0)[lb]{{\SetFigFont{8}{9.6}{\rmdefault}{\mddefault}{\updefault}$l'$}}}
\end{picture}
}
\caption{Midpoint}\label{f:midpoint}
\end{figure}

The {\em generic addition} configuration
(Figure~\ref{f:genericaddition})
deals with addition $x_a + x_b$ 
in the ``generic'' case where $x_a$, $x_b$, and $x_a+x_b$ are distinct
from $0$ and $1$, and the framing-type of their lines are all $\{ 0, 1 \}$.
 Given two lines $l_a$ and
$l_b$ with variables $x_a$, $x_b$, with framing points $(P_{a,0},
P_{a,1})$ and $(P_{b,0}, P_{b,1})$  on $l_a$
and $l_b$ respectively,  we choose a general horizontal line $l'$ and
a general point $X$, and superimpose the construction shown in Figure~\ref{f:genericaddition}.
If the line $l'$ is given  the framing-type $Fr = \{ 0,1 \}$ with framing
points $P'_0$ and $P'_1$, the reader will readily verify that the
coordinate of $V'$ is $x' = x_a+x_b$, and that this equation is
precisely what is (scheme-theoretically) enforced by the configuration.

\begin{figure}[ht]
\setlength{\unitlength}{0.00083333in}
\begingroup\makeatletter\ifx\SetFigFont\undefined%
\gdef\SetFigFont#1#2#3#4#5{%
  \reset@font\fontsize{#1}{#2pt}%
  \fontfamily{#3}\fontseries{#4}\fontshape{#5}%
  \selectfont}%
\fi\endgroup%
{\renewcommand{\dashlinestretch}{30}
\begin{picture}(3979,2011)(0,-10)
\put(1812,937){\blacken\ellipse{60}{60}}
\put(1812,937){\ellipse{60}{60}}
\put(3012,937){\blacken\ellipse{60}{60}}
\put(3012,937){\ellipse{60}{60}}
\put(1212,937){\blacken\ellipse{60}{60}}
\put(1212,937){\ellipse{60}{60}}
\put(1512,487){\blacken\ellipse{60}{60}}
\put(1512,487){\ellipse{60}{60}}
\put(2412,937){\blacken\ellipse{60}{60}}
\put(2412,937){\ellipse{60}{60}}
\put(612,37){\blacken\ellipse{60}{60}}
\put(612,37){\ellipse{60}{60}}
\put(1212,487){\blacken\ellipse{60}{60}}
\put(1212,487){\ellipse{60}{60}}
\put(1812,1387){\blacken\ellipse{60}{60}}
\put(1812,1387){\ellipse{60}{60}}
\put(2112,1387){\blacken\ellipse{60}{60}}
\put(2112,1387){\ellipse{60}{60}}
\put(2712,1387){\blacken\ellipse{60}{60}}
\put(2712,1387){\ellipse{60}{60}}
\put(2412,1837){\blacken\ellipse{60}{60}}
\put(2412,1837){\ellipse{60}{60}}
\put(1212,37){\blacken\ellipse{60}{60}}
\put(1212,37){\ellipse{60}{60}}
\put(2112,487){\blacken\ellipse{60}{60}}
\put(2112,487){\ellipse{60}{60}}
\dashline{60.000}(3612,1387)(12,1387)
\dashline{60.000}(3612,487)(12,487)
\path(3612,937)(12,937)
\path(3612,37)(12,37)
\path(612,37)(2412,937)
\path(1212,37)(1212,937)
\path(1212,37)(3012,937)
\path(612,37)(1812,937)
\path(2338,885)(2413,1035)(2488,885)(2338,885)
\path(2412,1837)(1212,937)
\path(1812,937)(2412,1837)
\path(2412,1837)(3012,937)
\path(2337,1912)(2487,1912)(2487,1762)
	(2337,1762)(2337,1912)
\path(3387,1012)(3537,1012)(3537,862)
	(3387,862)(3387,1012)
\put(3687,1387){\makebox(0,0)[lb]{{\SetFigFont{8}{9.6}{\rmdefault}{\mddefault}{\updefault}$l_a$}}}
\put(3687,937){\makebox(0,0)[lb]{{\SetFigFont{8}{9.6}{\rmdefault}{\mddefault}{\updefault}$l'$}}}
\put(3687,487){\makebox(0,0)[lb]{{\SetFigFont{8}{9.6}{\rmdefault}{\mddefault}{\updefault}$l_b$}}}
\put(1587,1462){\makebox(0,0)[lb]{{\SetFigFont{8}{9.6}{\rmdefault}{\mddefault}{\updefault}$P_{a,0}$}}}
\put(2112,1462){\makebox(0,0)[lb]{{\SetFigFont{8}{9.6}{\rmdefault}{\mddefault}{\updefault}$V_a$}}}
\put(2712,1462){\makebox(0,0)[lb]{{\SetFigFont{8}{9.6}{\rmdefault}{\mddefault}{\updefault}$P_{a,1}$}}}
\put(1062,1012){\makebox(0,0)[lb]{{\SetFigFont{8}{9.6}{\rmdefault}{\mddefault}{\updefault}$P'_0$}}}
\put(1662,1012){\makebox(0,0)[lb]{{\SetFigFont{8}{9.6}{\rmdefault}{\mddefault}{\updefault}$V$}}}
\put(2262,1012){\makebox(0,0)[lb]{{\SetFigFont{8}{9.6}{\rmdefault}{\mddefault}{\updefault}$V'$}}}
\put(912,562){\makebox(0,0)[lb]{{\SetFigFont{8}{9.6}{\rmdefault}{\mddefault}{\updefault}$P_{b,0}$}}}
\put(3012,1012){\makebox(0,0)[lb]{{\SetFigFont{8}{9.6}{\rmdefault}{\mddefault}{\updefault}$P'_1$}}}
\put(2187,1912){\makebox(0,0)[lb]{{\SetFigFont{8}{9.6}{\rmdefault}{\mddefault}{\updefault}$X$}}}
\put(2112,337){\makebox(0,0)[lb]{{\SetFigFont{8}{9.6}{\rmdefault}{\mddefault}{\updefault}$P_{b,1}$}}}
\put(1512,337){\makebox(0,0)[lb]{{\SetFigFont{8}{9.6}{\rmdefault}{\mddefault}{\updefault}$V_b$}}}
\end{picture}
}
\caption{Generic addition}\label{f:genericaddition}
\end{figure}
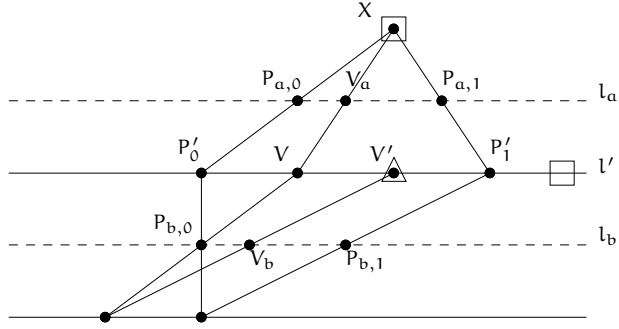

The {\em generic multiplication} configuration 
(Figure~\ref{f:genericmultiplication})
constructs/enforces
 multiplication $x_c = x_a  x_b$ in the ``generic''  case 
where $x_a$,  $x_b$, and $x_a x_b$ are distinct from $0$ and $1$, and
the framing-type of their lines are all $\{ 0, 1 \}$. 
As with the ``generic addition'' case, given two lines $l_a$ and $l_b$ with variables
$x_a$, $x_b$, with framing points $(P_{a,0},
P_{a,1})$ and $(P_{b,0}, P_{b,1})$  on $l_a$
and $l_b$ respectively,  we choose a general horizontal line $l'$ and
a general point $X$, and superimpose the construction shown in
Figure~\ref{f:genericmultiplication}.   
If the line $l'$ is given  the framing-type $Fr = \{ 0,1 \}$ with framing
points $P'_0$ and $P'_1$, the reader will readily verify that the
coordinate $x'$ of $V'$ is $x_a x_b$, and that this equation is
precisely what is (scheme-theoretically) enforced by the configuration.
The main part of the argument is that
\[
x_a = (P_{a,0}, V_a; P_{a,1}, p_3) = (P_0', V; P_1', p_3)
\]
and
\[
(P_0', V'; V, p_3) = (P_{b,0}, V_b; P_{b,1}, p_3) = x_b
\]
yield
\[
x' = (P_0', V'; P_1', p_3) = (P_0', V'; V, p_3) (P_0', V; P_1', p_3)  = x_a x_b.
\]
We remark that we are
parallel-shifting the point $V_b$ from $l_b$ to $l'$ to avoid
accidental overlaps of points in our later argument.

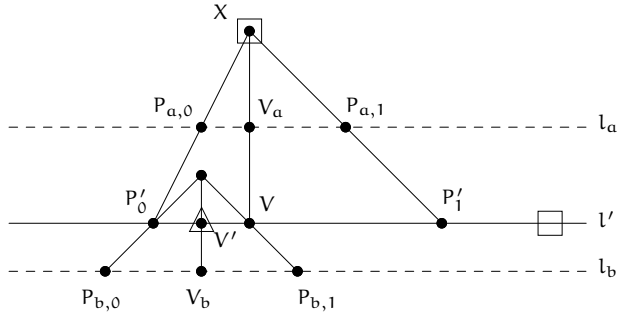
\begin{figure}[ht]
\setlength{\unitlength}{0.00083333in}
\begingroup\makeatletter\ifx\SetFigFont\undefined%
\gdef\SetFigFont#1#2#3#4#5{%
  \reset@font\fontsize{#1}{#2pt}%
  \fontfamily{#3}\fontseries{#4}\fontshape{#5}%
  \selectfont}%
\fi\endgroup%
{\renewcommand{\dashlinestretch}{30}
\begin{picture}(3979,1934)(0,-10)
\put(1212,860){\blacken\ellipse{60}{60}}
\put(1212,860){\ellipse{60}{60}}
\put(912,560){\blacken\ellipse{60}{60}}
\put(912,560){\ellipse{60}{60}}
\put(612,260){\blacken\ellipse{60}{60}}
\put(612,260){\ellipse{60}{60}}
\put(1212,560){\blacken\ellipse{60}{60}}
\put(1212,560){\ellipse{60}{60}}
\put(1212,260){\blacken\ellipse{60}{60}}
\put(1212,260){\ellipse{60}{60}}
\put(1512,560){\blacken\ellipse{60}{60}}
\put(1512,560){\ellipse{60}{60}}
\put(1812,260){\blacken\ellipse{60}{60}}
\put(1812,260){\ellipse{60}{60}}
\put(1212,1160){\blacken\ellipse{60}{60}}
\put(1212,1160){\ellipse{60}{60}}
\put(1512,1760){\blacken\ellipse{60}{60}}
\put(1512,1760){\ellipse{60}{60}}
\put(1512,1160){\blacken\ellipse{60}{60}}
\put(1512,1160){\ellipse{60}{60}}
\put(2112,1160){\blacken\ellipse{60}{60}}
\put(2112,1160){\ellipse{60}{60}}
\put(2712,560){\blacken\ellipse{60}{60}}
\put(2712,560){\ellipse{60}{60}}
\dashline{60.000}(12,1160)(3612,1160)
\path(12,560)(3612,560)
\dashline{60.000}(12,260)(3612,260)
\path(1437,1835)(1587,1835)(1587,1685)
	(1437,1685)(1437,1835)
\path(3312,635)(3462,635)(3462,485)
	(3312,485)(3312,635)
\path(1512,1760)(912,560)
\path(1512,1760)(1512,560)
\path(1512,1760)(2712,560)
\path(1212,860)(612,260)
\path(1212,860)(1212,260)
\path(1212,860)(1812,260)
\path(1138,509)(1213,659)(1288,509)(1138,509)
\put(462,35){\makebox(0,0)[lb]{{\SetFigFont{8}{9.6}{\rmdefault}{\mddefault}{\updefault}$P_{b,0}$}}}
\put(1812,35){\makebox(0,0)[lb]{{\SetFigFont{8}{9.6}{\rmdefault}{\mddefault}{\updefault}$P_{b,1}$}}}
\put(3687,1132){\makebox(0,0)[lb]{{\SetFigFont{8}{9.6}{\rmdefault}{\mddefault}{\updefault}$l_a$}}}
\put(3687,536){\makebox(0,0)[lb]{{\SetFigFont{8}{9.6}{\rmdefault}{\mddefault}{\updefault}$l'$}}}
\put(3687,227){\makebox(0,0)[lb]{{\SetFigFont{8}{9.6}{\rmdefault}{\mddefault}{\updefault}$l_b$}}}
\put(1114,31){\makebox(0,0)[lb]{{\SetFigFont{8}{9.6}{\rmdefault}{\mddefault}{\updefault}$V_b$}}}
\put(1573,638){\makebox(0,0)[lb]{{\SetFigFont{8}{9.6}{\rmdefault}{\mddefault}{\updefault}$V$}}}
\put(735,625){\makebox(0,0)[lb]{{\SetFigFont{8}{9.6}{\rmdefault}{\mddefault}{\updefault}$P'_0$}}}
\put(912,1230){\makebox(0,0)[lb]{{\SetFigFont{8}{9.6}{\rmdefault}{\mddefault}{\updefault}$P_{a,0}$}}}
\put(1563,1225){\makebox(0,0)[lb]{{\SetFigFont{8}{9.6}{\rmdefault}{\mddefault}{\updefault}$V_a$}}}
\put(2113,1234){\makebox(0,0)[lb]{{\SetFigFont{8}{9.6}{\rmdefault}{\mddefault}{\updefault}$P_{a,1}$}}}
\put(2712,635){\makebox(0,0)[lb]{{\SetFigFont{8}{9.6}{\rmdefault}{\mddefault}{\updefault}$P'_1$}}}
\put(1287,410){\makebox(0,0)[lb]{{\SetFigFont{8}{9.6}{\rmdefault}{\mddefault}{\updefault}$V'$}}}
\put(1287,1835){\makebox(0,0)[lb]{{\SetFigFont{8}{9.6}{\rmdefault}{\mddefault}{\updefault}$X$}}}
\end{picture}
}
\caption{Generic multiplication}\label{f:genericmultiplication}
\end{figure}

\section{Putting everything together}

We now put the atomic configurations together in various ways in order
to solve the four problems of \S \ref{s:four}.  We begin with the case
$p \neq 2$, leaving the case $p=2$ until \S \ref{s:char2}.

\bpoint{Relabeling} Before we start, we note that it will be convenient to use the same
framing points but a different framing-type to change the value of the
variable ``carried'' by the line. For example, a $(\{ 0,1 \},
q)$-line may be interpreted as a $( \{0, -1 \}, -q)$-line (as
$(0,1;q,\infty) = (0,-1; -q, \infty)$).\label{s:relabel}

\bpoint{Initial framing}
Before we start, we ``construct $-1$ on the $x$-axis''.  More precisely, on the $x$-axis, we have identified the points $0 :=
[0,0,1] = 
p_1$ and $1 := [1,0,1] = \overline{p_2 p_4} \cap \overline{p_1 p_3}$.
We use the midpoint construction (Figure~\ref{f:midpoint}) to construct
$-1 := [-1,0,1]$ as well (using $M=0$, $B=1$, $A=-1$).\label{s:if}

We now construct equality, negation, addition, and multiplication.

\bpoint{Equality:  enforcing $x_a=x_b$} We enforce equality $x_a=x_b$ as
follows.  \label{s:equality}

\epoint{First case: same framing-type} \label{s:e1}Suppose first that
two variables $x_a$ and $x_b$ are of same framing-type $\{ s_1, s_2
\}$. Then after a general choice of horizontal line $l'$, we parallel
shift (Figure~\ref{f:parallelshift}) $(P_{a, s_1}, V_a, P_{a, s_2})$
onto $(P_1', V', P_2')$ on $l'$, using a generally chosen projection point $X$. Then for $X' = P_1' P_{b, s_1} \cap P_2' P_{b, s_2}$, we shift $(P_{b, s_1}, V_b, P_{b, s_2})$ onto $(P'_1, V'',P_2')$ on $l'$, using $X'$ as projection point. The reader will verify that if we impose the
codimension $1$ condition that $V'=V''$, we enforce the equality $x_a
= x_b$.  The reader will verify that with the general choice of
projection point $X$ and line $l'$, the
newly constructed points will miss any finite number of previously
constructed points and lines (except for those in the Figure); and the
newly constructed lines will miss any finite number of previously
constructed points (except for those in the Figure, and of course
$p_3$) --- we will have no ``unintended coincidences''.  This can be
readily checked in all later constructions (an essential point in the
entire strategy!), but for concision's sake we will not constantly
repeat this.

\epoint{Second case:  different framing-type} \label{s:e2}Next, suppose that $x_a$ and $x_b$ have different framing-type, say
$\{ s_{a,1}, s_{a,2} \}$ and $\{ s_{b,1}, s_{b,2} \}$ respectively
(two distinct subsets of $\{ -1, 0, 1 \}$).  Then $q_a=q_b$ is not in $\{
-1, 0, 1 \}$.  We apply parallel shift (Figure~\ref{f:parallelshift})
to move $x_a$ to a generally chosen horizontal line $l'$.  We then
parallel shift the points $-1$, $0$, and $1$ on the $x$-axis to $l'$,
so we have marked points on $l'$ that can be identified (with the
obvious isomorphism to $\proj^1$) with $\infty=p_3$, $-1$, $0$, $1$,
and $q_a=q_b$.  Then (using the subset $\{ s_{b,1}, s_{b,2} \}$ of the
marked points on $l'$) $l'$ and $l_b$ have same framing-type and we can
apply previous construction.

We remark that in this and later constructions, we can take $x_a$ or
$x_b$ (or, later, $x_c$) to be the constants $0$ or $1$, by treating
the $x$-axis as a variable-bearing line.  For example, to take
$x_a\equiv 1$, treat the $x$-axis as a $( \{-1, 0 \}, 1)$-line.

\epoint{Remark:  choosing framing-type freely} The argument of \S \ref{s:e2}  shows that given a variable ``carried by'' a
variable-bearing line, we can change the framing-type of the line it
``lives on'', at the cost of moving it to another generally chosen horizontal line (so
long as the value of the variable does not lie in the new framing-type
of course).  From now on, given a variable, we freely choose
a framing-type to suit our purposes at the time.\label{r:e3}

\bpoint{Negation: enforcing $x_b=-x_a$} We now explain how to enforce $x_b =
-x_a$.  Suppose $x_a$ is carried on a line with framing-type $\{ s_1, s_2
\}$, and $x_b$ is carried on a line with framing-type $\{ -s_1, -s_2 \}$
(possible as $q_a=-q_b$ --- here we use Remark~\ref{r:e3}).  We
enforce $x_b=-x_a$ by adding the equality configuration (first case,
\S \ref{s:e1}), except interpreting line $l_b$ as an $( \{s_1, s_2, \},
-q_b)$-line (the relabel construction, \S \ref{s:relabel}).\label{s:n}

\bpoint{Addition:  enforcing $x_a+x_b=x_c$}\label{s:biga}

\epoint{First case: ``\text{(general)} + \text{(general)} = \text{(general)}''}\label{s:addition1}
Suppose   $q_a, q_b, q_c$ are all distinct from $0$ and $1$.
We apply  parallel shifts to move the three relevant variable-bearing
lines onto generally chosen lines, and then superimpose the ``generic
addition'' configuration of Figure~\ref{f:genericaddition}.
(The parallel shifts are to guarantee no unintended coincidences.)

\epoint{Second case: ``$1 + \text{(general)} = \text{(general)}$''}
\label{s:addition2}Suppose next that $q_a=1$, and $q_b$ and $q_b$ are
neither $0$ nor $1$.  Then $q_c \neq -1$ (or else $q_b$ would be $0$).
The equation we wish to enforce may be rewritten as $-x_c + x_b = -
x_a$, and $-x_c$, $x_b$, and $- x_a$ are all distinct from $0$ and
$1$. (Here we use $p \neq 2$, as we require $-1 \neq 1$.)  We thus accomplish our goal by applying the
negation configuration to $x_a$ and $x_c$, then applying 
the first case of the 
addition
construction, \S \ref{s:addition1}.

\epoint{Third case: ``$0 + \text{(general)} = \text{(general)}$''}
Suppose that $q_a=0$, and $q_b$ and $q_c$ are not in $\{ -1, 0, 1, 2
\}$.\label{s:addition3}
We take the framing-sets on $l_a$ and $l_c$ to be $\{ -1, 1 \}$ (using
Remark~\ref{r:e3}).
As in \S \ref{s:relabel}, we interpret/relabel
the $(\{ -1, 1 \}, x_a)$-line $l_a$ as a $(\{ 0, 1\} , x'_a)$-line
(where $x'_a = (x_a+1)/2$) and the $( \{-1, 1 \}, x_c)$-line
as a $(\{0, 1 \}, x'_c)$-line (where $x'_c = (x_c+1)/2$). 
We take the framing-set $\{0, 1 \}$ on $l_b$.
We parallel shift  $x_b$ onto a general horizontal line  $l'_b$, then use the midpoint 
construction (Figure~\ref{f:midpoint}) to construct the midpoint of
$V_b$ and $P_{b,0}$ on $l'_b$, so we have constructed the variable
$x_b/2$, which we name $x'_b$.  
The equation we wish to enforce, $x_a + x_b = x_c$, is algebraically equivalent to
$x'_a + x'_b = x'_c$, and the values of $x'_a$, $x'_b$, and $x'_c$ are
all distinct from $0$ and $1$, so we can apply the construction of the
first case of addition, \S \ref{s:addition1}.\label{s:a3}

\epoint{Fourth case: everything else}\label{s:addition4} We begin by
adding two extra free variables $s$ and $t$ on two generally chosen
horizontal lines.  More precisely, for $s$, we pick a generally chosen
horizontal line $l_i$, and three generally chosen points $P_{i,0}$,
$P_{i,1}$, and $V_i$ on it, 
and define  $s = (P_{i,0}, P_{i,1}; V_i, p_3)$,
so $l_i$ is a $(\{0,1 \}, s)$-line.  We do the same for $t$.  Using
the previous cases of addition, we successively construct $x_a + s$, $x_b + t$,
$(x_a + s) + (x_b + t)$, $s + t$, and $x_c + (s+t)$.  (Because $s$ and
$t$ were generally chosen, one of the three previous cases can always
be used.)  Then we impose
the equation
\[
(x_a + s) + (x_b + t) = x_c + (s+t)
\]
(using the third case of addition, \S \ref{s:a3}, twice).  
Thus we have scheme-theoretically enforced $x_a + x_b = x_c$ as
desired.

\bpoint{Multiplication: enforcing $x_a x_b = x_c$} \label{s:bigm} As
with addition, we deal with a ``sufficiently general'' case first, and
then deal with arbitrary cases by translating by a general value.

\epoint{First case:  ``$\text{(general)} \times
\text{(general)} = 
\text{(general)}$''}\label{s:multiplication1}
Suppose  $q_a, q_b, q_c \neq 0,1$.
We parallel shift all variables $x_a$, $x_b$, $x_c$ to generally chosen lines
$l'_a$, $l'_b$, and $l'_c$ (to avoid later unintended
incidences), and then  
superimpose the
generic multiplication configuration to impose $x_c=x_a x_b$
(where $l'_a$, $l'_b$, and $l'_c$ here correspond to $l_a$, $l_b$, and
$l'$ in Figure~\ref{f:genericmultiplication}).

\epoint{Second  case: everything else}
\label{s:multiplication2}To enforce $x_a x_b = x_c$, we  proceed as follows.
We add two extra free variables $u$ and $v$ as in \S
\ref{s:addition4}.  We then use the  addition constructions of \S \ref{s:addition1}--\ref{s:addition4} to construct $x_a + u$ and  $x_b + v$ (on generally
chosen 
horizontal lines).  We use the construction
of \S \ref{s:multiplication1} to construct $(x_a + u)(x_b + v)$, $uv$,
$(x_a + u)v$, and $(x_b+v) u$ (each on generally chosen lines).
Finally, we use the addition constructions (several times) to enforce 
$$(x_a + u)(x_b+v) + uv = x_c + (x_a+u)v + (x_b+v)u.$$
The result then follows from the algebraic identity
$$
(a+c) (b+d)  +cd = ab + (a+c)d + (b+d)c.
$$

\bpoint{Characteristic $2$}

\label{s:char2}As the above constructions at several points use $-1
\neq  1$, the case $p=2$ requires a variant strategy.

\epoint{Addition and multiplication: general cases (\S
  \ref{s:addition1}, \S \ref{s:multiplication1})} We begin by noting
that the general cases of addition and multiplication, given in \S
\ref{s:addition1} and \S \ref{s:multiplication1} respectively, work as
before (where $q_a$, $q_b$, and $q_c$ are all distinct from $\{ 0, 1
\}$, and the framing-type is taken to be $\{ 0, 1 \}$ in all cases).\label{s:2multiplication1}

\epoint{Relabeling (\S \ref{s:relabel}), and the first case of
  equality (\S \ref{s:e1})} Relabeling (\S \ref{s:relabel}) works as
before.  Equality in the case of same framing-type (\S \ref{s:e1})
does as well.  \label{s:2relabel}

\epoint{Initial framing} 
\label{s:if2}In analogy with the initial
framing of \S \ref{s:if},
before we begin the construction, we construct $j$ and $1-j=j^2$ on the
$x$-axis as follows.  More precisely, we will add whose points on the
$x$-axis which we label $j$ and $k$, as well as configurations forcing
the coordinates to satisfy $j^2+j-1=0$, and $k=j^2$.  (We then
hereafter call the point $k$ by the name $j^2$.)  It is important to note that
this construction of $j$ is \'etale over $\Spec \Z$ away from $[(5)]$,
and in particular at $2$; thus this choice will not affect the
singularity type.

We construct these points as follows.  Choose $j \in \F_4 \setminus
\F_2$, and place a marked point at $j$ on the $x$-axis.
Construct the product of $j$ with $j$ by parallel shifting $j$ separately onto two
generally chosen horizontal lines, and then using the construction of
\S \ref{s:2multiplication1}, i.e. \S \ref{s:multiplication1}  (possible as $j$ and $j^2$ are distinct from $0$ and
$1$).  Then construct $1-j$ using the relabeling trick of
\S \ref{s:relabel} (\S \ref{s:2relabel}): parallel shift $0$, $1$, and $j$ to a generally
chosen line, then reinterpret the $(\{0,1 \}, j)$-line as a $(\{ 1, 0
\}, 1-j)$-line, and parallel-shift it back to the $x$-axis.  Finally,
we use the equality configuration (the ``same framing-type'' case, \S
\ref{s:2relabel} = \S \ref{s:e1}) to
enforce $j^2 = 1-j$.

\epoint{Equality in general (\S \ref{s:e2}), and freely choosing
  framing-type (Remark~\ref{r:e3})}  Now that we have constructed $j$,
the second case of
the
equality construction works (with $\{ -1, 0, 1 \}$ replaced by $\{ 0,
1, j \}$), and we may choose framing-type freely on lines as observed
in Remark~\ref{r:e3}.

\epoint{Addition, second case: ``$1 + \text{(general)} =
  \text{(general)}$'', cf.\  \S \ref{s:addition2}}
Suppose $q_a=1$, and $q_b$ and $q_c$ are not in $\{0, 1, j \}$.
Then use  the general case of multiplication (\S
\ref{s:2multiplication1}, i.e.\ \S \ref{s:multiplication1}) to construct (on separate general horizontal
lines)
$q'_b = q_b/j$ and $q'_c = q_c/j$.  By considering the $( \{ 0, j \},
x_a)$-line
as a $( \{ 0,1 \}, x_a/j)$-line (\S \ref{s:2relabel}, i.e.\ \S \ref{s:relabel}), construct (on
a  general  horizontal line, using parallel shift) $x'_a = x_a/j$.  Then impose $x'_a + x'_b
= x'_c$ using the general case of addition (\S
\ref{s:2multiplication1})

\epoint{Addition, third case: ``$0 + \text{(general)} =
  \text{(general)}$'', cf.\ \S \ref{s:addition3}} Suppose $q_a = 0$,
$q_b \notin \{ 0, 1, j, j^2 \}$, and $q_c \notin \{ 1, j \}$.
Then construct  $x'_a = (x_a-1)/j^2$
(on a general horizontal line of framing-type $\{0, 1 \}$)
 by
considering the $( \{ 1, j \}, x_a)$-line (carrying the variable
$x_a$) as a $( \{ 0, 1 \}, (x_a-1)/j^2)$-line (as $j-1=j^2$).
Similarly, construct $x'_c = (x_c-1)/j^2$.  Using the general
multiplication construction (\S \ref{s:2multiplication1}, i.e.\ \S \ref{s:multiplication1}) twice, construct
$x'_b = x_b/j^2$ (by way of the intermediate value of $x_b /j$).  Then
impose $x'_a + x'_b = x'_c$ (using the general addition construction
of \S \ref{s:2multiplication1}, i.e.\ \S \ref{s:addition1}), and note that this is algebraically
equivalent to $x_a + x_b = x_c$.

\epoint{Addition and multiplication, final cases:  everything else (\S
  \ref{s:addition4} and \S \ref{s:multiplication2})}
These now work as before.

\epoint{Negation (\S \ref{s:n})} Finally, negation can be imposed by
constructing the configuration imposing $x_a + x_b = 0$ (using the
final case of addition).



} 


\begin{thebibliography}{[XX]}
\bibitem[BB]{bb} P. Belkale and P. Brosnan, {\em Matroids, motives,
    and a conjecture of Kontsevich}, Duke Math.\ J.\ {\bf 116} (2003),
  no.\  147--188.
\bibitem[BS]{bs} J. Bokowski and B. Sturmfels, {\em Computational
    Synthetic Geometry}, Lecture Notes in Math.\ {\bf 1355}, Springer,
  Berlin, 1989.
\bibitem[E]{e} D. Erman, {\em Murphy's Law for Hilbert function strata
    in the Hilbert scheme of points}, arXiv:1205.0587.
\bibitem[HM]{hm}  J. Harris and I. Morrison, {\em Moduli of Curves},
  GTM {\bf 187}, Springer-Verlag, New York, 1998.
\bibitem[KT]{kt} S. Keel and J. Tevelev, {\em Geometry of Chow
    quotients of Grassmannians}, Duke Math.\  J.\ {\bf 134} (2006),
  no.\ 2, 259--311.
\bibitem[Ku]{k}  J. P. S. Kung, {\em A Source Book in Matroid Theory},
  Birkh\"auser, Boston, 1986.
\bibitem[L]{l} L. Lafforgue, {\em Chirurgie des Grassmanniennes}, CRM
  Monograph Series {\bf 19}, Amer.\ Math.\ Soc., 2003.
\bibitem[Ma]{mac} S. Mac Lane, {\em Some interpretations of abstract
    linear dependence in terms of projective geometry}, Amer.\ J.\
  Math.\ {\bf 58} (1936), 236--240.
\bibitem[M1]{m1} N. Mn\"ev, {\em Varieties of combinatorial types of
    projective configurations and convex polyhedra}, Dokl.\ Akad.\
  Nauk SSSR {\bf 283}, 1312--1314 (1985).
\bibitem[M2]{m2} N. Mn\"ev, {\em The universality theorems on the
    classification problem of configuration varieties and convex
    polytopes varieties}, in {\em Topology and Geometry --- Rohlin
    Seminar}, LNM {\bf 1346}, 527--543, Springer,
  Berlin, 1988.
\bibitem[P]{p} S. Payne, {\em Moduli of toric vector bundles},
  Compos.\ Math.\ {\bf 144} (2008), no.\ 5, 1199--1213.
\bibitem[S1]{s1} B. Sturmfels, {\em On the decidability of Diophantine
    problems in combinatorial geometry}, Bull.\ Amer.\ Math.\ Soc.\
  (N.S.) {\bf 17} (1987), no.\ 1, 121--124.
\bibitem[S2]{s2} B. Sturmfels, {\em On the matroid stratification of
    Grassmann varieties, specialization of coordinates, and a problem
    of N. White}, Adv.\ Math.\ {\bf 75} (1989), 202--211.
\bibitem[Va]{v} R. Vakil, {\em Murphy's law in algebraic geometry:
    Badly-behaved deformation spaces}, Invent.\ Math.\ {\bf 164} (2006),
  569--590.
\bibitem[Ve]{vershik} A. M. Vershik, {\em Topology of the convex
    polytopes' manifolds, the manifold of the projective
    configurations of a given combinatorial type and representations
    of lattices}, in {\em Topology and Geometry --- Rohlin
    Seminar}, LNM {\bf 1346}, 557--581, Springer,
  Berlin, 1988.
\end{thebibliography}
\end{document}